\documentclass[reqno]{amsart}
\usepackage{amsmath,amssymb,cite}
\usepackage{graphics,epsfig,color}
\usepackage[mathscr]{euscript}
\usepackage{mathtools}
\mathtoolsset{showonlyrefs}
\usepackage{xcolor}
\usepackage{braket}

\usepackage{ulem}

\definecolor{bblue}{rgb}{.2,0.2,.8}

\theoremstyle{plain}
\newtheorem{theorem}{Theorem}[section]

\theoremstyle{definition}
\newtheorem{definition}[theorem]{Definition}

\theoremstyle{remark}
\newtheorem{remark}[theorem]{Remark}

\numberwithin{equation}{section}
\numberwithin{theorem}{section}

\def\be{\begin{equation}}
\def\ee{\end{equation}}
\def\bp{\begin{pmatrix}}
\def\ep{\end{pmatrix}}
\def\bea{\begin{eqnarray}}
\def\eea{\end{eqnarray}}

\def\\{\par\medskip}

\let\0=\noindent

\newcommand{\id}{{1 \mskip -5mu {\rm I}}}
\renewcommand{\epsilon}{\varepsilon}
\renewcommand{\hat}{\widehat}

\begin{document}

\title[Mixtures and the matrix product ansatz]{Mixtures, Markov bridges and the matrix product ansatz}

\author{Davide Gabrielli}
\address{\noindent Davide Gabrielli \hfill\break\indent
 DISIM, Universit\`a dell'Aquila
\hfill\break\indent
67100 Coppito, L'Aquila, Italy
}
\email{davide.gabrielli@univaq.it}

\author{Federica Iacovissi}
\address{\noindent Federica Iacovissi \hfill\break\indent
	DISIM, Universit\`a dell'Aquila
	\hfill\break\indent
	67100 Coppito, L'Aquila, Italy
}
\email{federica.iacovissi@graduate.univaq.it}

\begin{abstract}
We give a probabilistic characterization of the set of measures that can be represented by the matrix product ansatz. By suitably enlarging the state space, we show that a probability measure can be described in terms of non negative matrices by the {\it Matrix Product Ansatz}, if and only if it can be written as a mixture of inhomogeneous product measures where the mixing law is a Markov bridge. We give a constructive procedure to identify such probabilistic features. We illustrate the result by examples and show that existing probabilistic representations of the invariant measures of non equilibrium interacting particle systems can be obtained from the matrix product ansatz by this general procedure.
\end{abstract}

\noindent


\maketitle
\thispagestyle{empty}

\section{Introduction}
Several important non equilibrium models of particles systems have an invariant measure that can be written in terms of product of matrices, this is the so called {\it Matrix Product Ansatz} (MPA). The MPA has a long history starting from \cite{Der1,Der2,Der3,Schutz}, with a large number of applications to many different models (see \cite{reviewhold,review} for some review of the method and \cite{D} for a review of some applications). The MPA is related to several combinatorial results and constructions, we here instead explore its probabilistic structure. We point out that measures of MPA type are studied also in informatics and pattern statistics under the name of rational models, see for example \cite{G2}.

\smallskip

Our main result is that a probability measure can be written by the MPA using matrices with non negative entries, if and only if the measure has a special canonical probabilistic structure. This special structure  is a mixture of non homogeneous product measures when the parameters that determine the marginals of the product are distributed according to a Markov bridge. A Markov bridge is a time inhomogeneous Markov chain that is informally obtained from an homogeneous one by a conditioning on the joint law of initial and final values.

The result is obtained enlarging the state space by a special coupling and then using a generalized Doob transform that associates a stochastic matrix to a non negative one. This transformation is based on a spectral construction. We prove our general statement in the case of finite alphabets, when the spectral structure is easier and we do not need to discuss special and detailed assumptions, but the basic mechanism holds also in the infinite setting. We illustrate indeed the construction by several examples of particle systems where the matrices involved in the MPA are infinite dimensional. In particular we show that some remarkable probabilistic representations of the invariant measures of non reversible particle systems, like the two lines ensemble for the boundary driven TASEP \cite{BC,Bryc}, or the queue representation of a multiclass TASEP \cite{A, FM, EFM}, can be obtained from the corresponding MPA by using our transformation.

Recently representations in terms of mixtures of inhomogeneous product measures of the invariant measures of boundary driven generalized KMP and zero range models have been obtained in \cite{CFGGT,dMFG} and in \cite{CFF, GRV,RV}. For the zero range case, considering just the simplest model, we show (it is one of the examples we consider) that this representation is a special case of our construction corresponding to a Markov process and marginals laws of the product measures having special features. While in the general case the marginals of the inhomogeneous product measure depend on two values of the Markovian parameter hidden variables, for these models the dependence reduces to just one. Moreover in this case the Markov bridges are associated to increasing processes so that the Markov property implies the spatial Markov property of \cite{RV}. In the case of the KMP and generalized KMP processes it is not clear if the hidden variables that determine the parameters of the inhomogeneous products are Markovian, so that is not clear if the invariant measures are of MPA type.

\smallskip

Our enlarging state space and a consequent Markov characterization are similar to the construction, in a stationary setting, of the algebraic and manifestly positive algebraic measures in \cite{FNS}. In \cite{FNS} a characterization of such measures it is obtained as function of Markov processes for a sufficiently enlarged state space, see our remark \ref{brno} for more details. Our mixture representation does not require an enlargement of the state space and it is obtained by a further step using the fact that the Markov processes involved are of special type. 

It would be also interesting to apply our construction to models on the infinite line \cite{SW}, models with finite dimensional representations \cite{MS, Sc} and to compare our construction with recent results for quantum models \cite{S}.

\smallskip
Our mixture representation is very suitable to handle the measures and in particular to deduce large deviations asymptotic. This will be done in \cite{GI} and \cite{I} where also other examples will be discussed.

\smallskip

The paper is organized as follows.

In section 2 we fix notation and introduce the basic combinatorial and probabilistic objects discussed, like the MPA, the Markov bridges, and the mixtures. At the end of the section we state the main Theorem of equivalence between MPA and mixtures of inhomogeneous product measures with parameters distributed as a Markov bridge. We state our Theorem in the simplest finite dimensional setting.

In section 3 we prove the main Theorem and discuss some remarks on the construction and related results.

In section 4 we discuss some examples; the majority of them are infinite dimensional.

\section{Notation and results}
We consider two finite or countable sets $A$ and $B$. The set $A$ is the state space where the variables that we study assume values while $B$ is an auxiliary state space. We discuss notation using $A$ but a similar notation holds for $B$ too and any other set.

Given $\eta=(\eta_i)_{i=1}^N, \in A^N$, we denote by $\eta_i^j=(\eta_i,\dots ,\eta_j)\in A^{j-i+1}$, the finite portion of the word $\eta$ contained between the indices $i<j$. Given $\eta,\gamma\in A^N$ we denote by $\eta_i^j\gamma_l^m=(\eta_i,\dots ,\eta_j,\gamma_l\dots ,\gamma_m)\in A^{j+m-i-l+2}$ the concatenation.  Given a set $S$ we denote by $\mathcal M^1(S)$ the probability measures on $S$. Given $\mu_N\in \mathcal M^1(A^N)$ and $\eta\in A^N$ we denote by $\mu_N(\eta_i^j):=\sum_{\{\gamma\in A^N:\gamma_i^j=\eta_i^j\}}\mu_N(\gamma)$ the probability of the cylinder set associated to $\eta_i^j$. The probability of the cylinder set determined by the single value $\eta_i$ is denoted by $\mu_N(\eta_i)$.

\subsection{The matrix product ansatz}

We introduce a family of probability measures, that we call of {\it Matrix Product Anstatz} (MPA) type, that form a special subset of $\mathcal M^1(A^N)$.
For any element $a\in A$ we have a $B\times B$ matrix $M^a$ having non-negative entries, i.e. for any $b,b'\in B$ we have $M^a_{b,b'}\geq 0$. We consider also $x:=(x_b)_{b\in B},y:=(y_b)_{b\in B}\in \mathbb R_+^B$ two column vectors having non negative entries, and denote by $x^T, y^T$ the corresponding row vectors obtained by transposition.
\begin{definition}
An element $\mu_N\in \mathcal M^1(A^N)$ is of {\it Matrix Product Anstatz} (MPA) type if it can be written as
\begin{equation}\label{defbase}
\mu_N(\eta):=\frac{y^T \left(\prod_{i=1}^NM^{\eta_i}\right)x}{Z_N}\,, \qquad \eta\in A^N\,,
\end{equation}
where $x,y$ and $(M^a)_{a\in A}$ are non negative vectors and matrices and
$Z_N$ is a normalization factor. The product of the matrices and vectors is the usual lines by column product.
\end{definition}
The probability measure $\mu_N$ depends on the family of matrices and the vectors but we will explicit such a dependence only when necessary. When $|A|,
|B|<+\infty$ the measure \eqref{defbase} is always well defined while in the countable infinite cases the measure exists only when $Z_N$ is finite.

We call $M:=\sum_{a\in A}M^a$ and the normalization term is given by
\begin{equation}\label{norm}
Z_N=\sum_\eta y^T \left(\prod_{i=1}^NM^{\eta_i}\right)x= y^T M^N x\,.
\end{equation}

\subsection{Markov bridges}

Consider a $B\times B$ stochastic matrix $P$ and two column vectors $f,g\in \mathbb R_+^B$. For notational convenience, given a $B\times B$ stochastic matrix $P$ and $\zeta\in B^{N+1}$, we denote by
$$\mathbb P_P(\zeta):=\prod_{i=1}^NP_{\zeta_i,\zeta_{i+1}}\,,$$
i.e. the probability of the path $\zeta$ with fixed initial condition $\zeta_1$ according to the Markov law induced by $P$ (a similar notation is also used for other state space). We define then a family of elements $\rho_{N+1}\in\mathcal M^1(B^{N+1})$ by
\begin{equation}\label{defbridge}
\rho_{N+1}(\zeta):= \frac{g_{\zeta_1}\mathbb P_P(\zeta)f_{\zeta_{N+1}}}{g^T P^N f}\,, \qquad \zeta\in B^{N+1}\,.
\end{equation}
The notation is like for the formula \eqref{defbase} and again we do not write explicitly the dependence of the measure on the parameters $f,g,P$. The motivation of the $N+1$ instead of $N$ will be clarified later on. The term in the denominator is the normalization factor. Again in the finite case the measure \eqref{defbridge} is always well defined while a summability condition is necessary in the infinite case.

\smallskip
The measure \eqref{defbridge} is a non-homogeneous Markov measure since by a direct computation we have
\begin{equation}\label{prtrbr}
\mathbb \rho_{N+1}(\zeta_{k+1}|\zeta_1^k)=\frac{P_{\zeta_k,\zeta_{k+1}}\sum_{\xi_{N+1}}P^{N-k}_{\zeta_{k+1},\xi_{N+1}}f_{\xi_{N+1}}}{\sum_{\xi_{N+1}}P^{N+1-k}_{\zeta_{k},\xi_{N+1}}f_{\xi_{N+1}}}\,, \qquad \zeta\in B^{N+1}\,.
\end{equation}
The Markovianity follows by the fact that the right hand side in the above formula does not depend on $\zeta_1^{k-1}$, i.e. $\rho_{N+1}(\zeta_{k+1}|\zeta_1^k)=\rho_{N+1}(\zeta_{k+1}|\zeta_k)$.

\smallskip

In the cases that $f_{\cdot}=\delta_{\cdot, b'}$ and $g_{\cdot}=\delta_{\cdot, b}$ the Markov bridge \eqref{defbridge} is pinned at $b$ at time $1$ and at $b'$ at time $N+1$. In this case \eqref{defbridge} is the law of the Markov chain with transition probability $P$ and conditioned to start at $b$ and to be in $b'$ at time $N+1$.
We call $\rho_{N+1}^{b,b'}$ this measure that is given by
\begin{equation}\label{defbr}
\rho_{N+1}^{b,b'}(\zeta)=\left\{
\begin{array}{ll}
\frac{\mathbb P_p(\zeta)}{P^N_{b,b'}}\,, & \textrm{if} \ \zeta_1=b\,, \zeta_{N+1}=b'\,,\\
0\,, & \textrm{otherwise}\,.
\end{array}
\right.
\end{equation}
The above measure depends of course also from the stochastic matrix $P$ but we do not write explicitly such dependence.

In the general case we have that the joint law of the pair $(\zeta_1,\zeta_{N+1})$ is given by $m\in \mathcal M^1(B^2)$ defined by
\begin{equation}\label{emme}
m(\zeta_1,\zeta_{{N+1}})=\frac{g_{\zeta_1} P^N_{\zeta_1,\zeta_{N+1}} f_{\zeta_{N+1}}}{g^T P^N f}
\end{equation}
and then we have the convex decomposition $$\rho_{N+1}(\zeta)=\sum_{b,b'}m(b,b')\rho_{N+1}^{b,b'}(\zeta)\,.$$
This means that a sample of the measure \eqref{defbridge} can be generated by choosing initial and final values according to $m$ in \eqref{emme}, and then generating a pinned Markov bridge with those initial and final values.

\subsection{Mixtures}
Consider a map $\hat p:\{1,\dots ,N\}\times B^{N+1}\to \mathcal M^1(A)$, that associates to the pair $i\in\{1,\dots , N+1\}$ and $\zeta\in B^{N+1}$ the probability measure $\hat p_\zeta^i(\cdot)\in \mathcal M^1(A)$.

\begin{definition}\label{defmix}
	We say that a probability measure $\nu_N\in \mathcal M^1(A^N)$ is a mixture of product measures with marginal measures determined by $\hat p$ and parameter distribution $\rho_{N+1}\in \mathcal M^1(B^{N+1})$ if we have
	\begin{equation}\label{mixform}
	\nu_N(\eta)=\sum_{\zeta\in B^{N+1}} \rho_{N+1}(\zeta)\left[\prod_{i=1}^N\hat p^i_\zeta(\eta_i)\right]\,, \qquad \eta\in A^N\,.
	\end{equation}
\end{definition}
We use the symbol $\rho_{N+1}$ for the parameter distribution that is the same symbol we use for a Markov bridge just for simplicity of notation. In Definition \ref{defmix} the measure $\rho_{N+1}\in \mathcal M^1(B^{N+1})$ is arbitrary, but in the following we will consider the special cases when it is indeed a Markov bridge.

In principle we could consider parameter distributions on longer sequences in the alphabet $B$ but the length $N+1$ is natural due to the fact that we will study the special cases when
\begin{equation}\label{spform}
\hat p^i_\zeta(\cdot)=p_{\zeta_i,\zeta_{i+1}}(\cdot)\,,
\end{equation}
for a suitable map $p: B^2\to \mathcal M^1(A)$ that associate to the pair $(b,b')\in B^2$ the probability measure $p_{b,b'}(\cdot)\in \mathcal M^1(A)$.

\smallskip
Comments on the definition \ref{defmix} and the related constructions are postponed to the following sections and we move directly to the main result.

\subsection{The equivalence}

We state and prove our main result in the cases of finite alphabets where we can avoid any technical assumption; we will then illustrate the result with examples that are not finite, since the basic construction works also in that case modulo some technical assumptions, that are satisfied in the examples we consider.

Recall by the Perron-Frobenius theorem (see for example \cite{K}) that a finite irreducible matrix $M$ having non negative entries has a positive maximal eigenvalue $\lambda$ and a corresponding strictly positive right eigenvector $e$ such that $Me=\lambda e$. By a simple transformation, sometimes called a generalized Doob transform or ground state transform, we associate to $M$ a stochastic matrix given by
\begin{equation}\label{Doob}
P=\lambda^{-1}\mathcal E^{-1}M\mathcal E\,,
\end{equation}
where $\mathcal E$ is the diagonal matrix having elements $\mathcal E_{i,i}=e_i$. The entries of the matrix $P$ are then given by $P_{b,b'}=\frac{M_{b,b'}e_{b'}}{\lambda e_b}$.

This construction is also related to Hammersley–Clifford theorem.

\begin{theorem}\label{ilteo}
	Consider $|A|,|B|<+\infty$; a probability measure $\mu_N \in \mathcal M^1(A^N)$ is of MPA type \eqref{defbase} with $M:=\sum_{a\in A}M^a$ an irreducible non negative matrix, if and only if $\mu_N$ is a mixture of product measures \ref{defmix} having the distribution of the parameters $\rho_{N+1}$ that is a Markov bridge \eqref{defbridge} and marginal measures of the form \eqref{spform}. Moreover the Markov bridge has transition probability given by the stochastic matrix $P$ defined in \eqref{Doob}; the vectors $f,g$ are related to $x,y$ (those of the MPA) by the relations $f=\mathcal E^{-1}x$, $g=\mathcal E y$ and the marginal distributions are of the form \eqref{spform} with
	\begin{equation}\label{defp}
	p_{b,b'}(a)=\frac{M^a_{b,b'}}{M_{b,b'}}\,, \qquad b,b'\in B\,,\, a\in A\,.
	\end{equation}
	Finally we have the relation $g^TP^Nf=\frac{Z_N}{\lambda^N}$.
\end{theorem}

\section{Proof of Theorem \ref{ilteo} and remarks}

The proof is based on a special coupling construction and the transform \eqref{Doob}. In the proof we will introduce a coupling measure and discuss some of its properties.

\begin{proof}{\it  [of Theorem \ref{ilteo}]}
	
First we start with a measure of MPA type \eqref{defbase} and show that is a mixture. Consider a measure $\mu_N$ like in \eqref{defbase}. We enlarge the state space and construct a probability measure
	$C_N\in \mathcal M^1(A^N\times B^{N+1})$ defined by
	\begin{equation}\label{coupling}
	C_N(\eta,\zeta):=\frac{y_{\zeta_1}\left(\prod_{i=1}^NM^{\eta_i}_{\zeta_i,\zeta_{i+1}}\right)x_{\zeta_{N+1}}}{Z_N}\,, \qquad \eta\in A^N\,,\, \zeta\in B^{N+1}\,.
	\end{equation}
	The coupling measure $C_N$ has the following remarkable properties.
	\begin{enumerate}
		\item The $\eta$ marginal of $C_N$ is $\mu_N$ in \eqref{defbase}, i.e. $\mu_N(\eta)=\sum_{\zeta\in B^{N+1}} C_N(\eta,\zeta)$.
		\item The measure $C_N$ is a Markov bridge with state space $A\times B$.
		\item The $\zeta$ marginal of $C_N$, $\rho_{N+1}(\zeta)=\sum_{\eta\in A^N} C_N(\zeta,\eta)$, is a Markov bridge with the features described in the text of the theorem.
		\item Conditioned on the $\zeta$ variables, the $\eta$ variables are independent and moreover
		\begin{equation}
		C_N(\eta|\zeta)=\prod_{i=1}^Np_{\zeta_i,\zeta_{i+1}}(\eta_i)\,,
		\end{equation}
		where the $p_{b,b'}(\cdot)$ are defined in \eqref{defp}.
	\end{enumerate}
	
	Item $(1)$ follows directly by definition since
	$$
	\sum_{\zeta\in B^{N+1}} y_{\zeta_1}\left(\prod_{i=1}^NM^{\eta_i}_{\zeta_i,\zeta_{i+1}}\right)x_{\zeta_{N+1}}=y \prod_{i=1}^NM^{\eta_i}x\,.
	$$
	
	\smallskip
	
	Item $(2)$ will not be used in this paper and we skip some details related to the irreducibility of the
	positive matrix we are going to introduce, these details will be explained in \cite{GI}. We introduce the $(A\times B)\times (A\times B)$ matrix $T_{(a,b),(a',b')}:=M^a_{b,b'}$ that has non negative entries. As will be shown in \cite{GI} if $M$ is irreducible then $T$ is irreducible too. By \eqref{Doob} we can construct the $(A\times B)\times (A\times B)$ stochastic matrix $S:=\Lambda^{-1}\hat{\mathcal  E}^{-1}T\hat {\mathcal E}$ where $\Lambda$ is the Perron eigenvalue of $T$ and ${\mathcal E}$ is the diagonal matrix having elements $\hat{\mathcal E}_{(a,b),(a,b)}=\varepsilon_{a,b}$, where $\varepsilon$ is the positive eigenvector associated to the eigenvalue $\Lambda$. Introduce the vectors $\hat x, \hat y\in \mathbb R_+^{A\times B}$ defined by $\hat x_{a,b}:=x_b$ and $\hat y_{a,b}.=y_b$. We can write then
	$$
	C_N(\eta,\zeta)=\frac{\hat y_{\eta_1,\zeta_1}\prod_{i=1}^NT_{(\eta_i,\zeta_i),(\eta_{i+1},\zeta_{i+1})}\hat x_{\eta_{N+1},\zeta_{N+1}}}{Z_N}\,,
	$$
	where we observe that in the above formula, differently from \eqref{defbase}, there is only one matrix appearing (that is $T$). We stress also that even if in the right hand side of the formula above appears the variable $\eta_{N+1}$, that is not defined and does not appear on the left hand side, due to the special form of $T$ and $\hat x$ also the right hand side is independent of $\eta_{N+1}$. By a telescoping argument we can write the above formula as
	$$
	C_N(\eta,\zeta)=\frac{\Lambda^N\hat y_{\eta_1,\zeta_1}\hat{\mathcal E}_{(\eta_1,\zeta_1),(\eta_1,\zeta_1)}\mathbb P_S(\eta,\zeta) \hat{\mathcal E}^{-1}_{(\eta_{N+1},\zeta_{N+1}),(\eta_{N+1},\zeta_{N+1})}\hat x_{\eta_{N+1},\zeta_{N+1}}}{Z_N}\,,
	$$
	that has exactly the form \eqref{defbridge} but on the state space $A\times B$ with the transition matrix $P$ given by $S$ and the vectors $g,f$ given by $\hat {\mathcal E}\hat y\,,\, \hat {\mathcal E}^{-1}\hat x$.
	
	\smallskip
	
	Item $(3)$ is obtained by
	\begin{equation}\label{lollo}
	\sum_\eta y_{\zeta_1}\left(\prod_{i=1}^NM^{\eta_i}_{\zeta_i,\zeta_{i+1}}\right)x_{\zeta_{N+1}}=y_{\zeta_1}\left(\prod_{i=1}^NM_{\zeta_i,\zeta_{i+1}}\right)x_{\zeta_{N+1}}\,.
	\end{equation}
	Using now formula \eqref{Doob}, by the same telescoping argument as before we obtain
	\begin{align}
	\rho_{N+1}(\zeta)&=\sum_\eta C_N(\eta,\zeta)\\
	&=\frac{y_{\zeta_1}\left(\prod_{i=1}^NM_{\zeta_i,\zeta_{i+1}}\right)x_{\zeta_{N+1}}}{Z_N}\\
	&=\frac{\lambda^N\left(\mathcal Ey\right)_{\zeta_1}\mathbb P_P(\zeta)\left(\mathcal E^{-1} x\right)_{\zeta_{N+1}}}{Z_N}\,,
	\end{align}
	that is of the form \eqref{defbridge} and we deduce the validity of the statement of item $(3)$.
	
	\smallskip
	
	Item $(4)$ is obtained by a direct computation using \eqref{coupling} and \eqref{lollo}
	$$
	C_N(\eta|\zeta)=\frac{C_N(\eta,\zeta)}{\rho_{N+1}(\zeta)}=\prod_{i=1}^N\frac{M^{\eta_i}_{\zeta_i,\zeta_{i+1}}}{M_{\zeta_i,\zeta_{i+1}}}=\prod_{i=1}^Np_{\zeta_i,\zeta_{i+1}}(\eta_i)\,.
	$$
	The relation $g^TP^Nf=\frac{Z_N}{\lambda^N}$ follows directly by \eqref{norm} and \eqref{Doob}.
	
	\medskip
	
	Conversely consider now a measure written in terms of mixtures like in Definition \ref{defmix} with the measure $\rho_{N+1}$ being a Markov bridge \eqref{defbridge} and the marginal distributions that are like in \eqref{spform}. We introduce the vectors $x=f$ and $y=g$ and the matrices $M^a$ whose elements are defined by $M^a_{b,b'}=P_{b,b'}p_{b,b'}(a)$. We obtain in this way a representation by the MPA of the mixture; with this choice the matrix $M$ is stochastic and coincides with $P$.
	
\end{proof}

\begin{remark}\label{Remo1}
The mixture representation \eqref{mixform}, without imposing some restrictions on the measures involved, is a very general statement and can be shown under general assumptions that such a representation of a probability measure $\mu_N$ exists and it is moreover not unique. Consider for example $A=B=\{0,1\}$ and the marginal distributions given by $\hat p^i_\zeta(\cdot)=\mathcal B_{\zeta_i}(\cdot)$ where $\mathcal B_p(\cdot)$ is the Bernoulli distribution of parameter $p$. We have for example that any $\mu_N\in \mathcal M^1\left(\{0,1\}^N\right)$ has always a representation like \eqref{mixform} as
$$
\mu_N(\eta)=\sum_{\zeta\in \{0,1\}^{N+1}}\rho_{N+1}(\zeta)\left[\prod_{i=1}^NB_{\zeta_i}(\eta_i)\right]\,,
$$
provided that $\rho_{N+1}(\zeta_1^N)=\mu_N(\zeta_1^N)$.

\smallskip

To illustrate the absence of uniqueness let us consider the following simple case.
Consider the product measure $\mu_N\in \mathcal M^1\left(\{0,1\}^N\right)$ given by $\mu_N=\prod_{i=1}^N\mathcal B_{q_i}$ for an arbitrary collection of parameters $q_i\in [0,1]$. Consider
$\rho_{N+1}\in \mathcal M^1\left(B^{N+1}\right)$ that is also a product measure with $i-$marginal given by $\rho^{[i]}\in \mathcal M^1(B)$. Consider functions $m_i:B\to[0,1]$ with mean values $\sum_{b\in B}\rho^{[i]}(b)m_i(b)=q_i$ and consider the marginal distributions $\hat p^i_\zeta(\cdot)=\mathcal B_{m_i(\zeta_i)}(\cdot)$. We have under this only assumption that
$$
\prod_{i=1}^N\mathcal B_{q_i}(\eta_i)=\sum_{\zeta\in B^{N+1}}\rho_{N+1}(\zeta)\left[\prod_{i=1}^N\mathcal B_{m_i(\zeta_i)}(\eta_i)\right]\,.
$$
Since the condition is only on the mean values of the marginals of $\rho_{N+1}$ it is easy to see that uniqueness does not hold.
\end{remark}

\begin{remark}
The importance and power of the representation as a mixture \eqref{mixform} for measures of MPA type is on the fact that the law $\rho_{N+1}$ of the parameters is Markovian and the marginal distributions $\hat p^i_\zeta(\cdot)$ are local (i.e. depends just on $\zeta_i,\zeta_{i+1}$) and translational covariant. In particular mean values, covariances and other properties of the measure $\mu_N$ are strictly related to those of the underlying Markov bridge. Consider for example correlations. Let us define $V_{\zeta_i,\zeta_{i+1}}:=\sum_{a\in A}p_{\zeta_i,\zeta_{i+1}}(a)a$. By a direct computation we have
$$
\mathbb E_{\mu_N}\left(\eta_i\eta_j\right)=\mathbb E_{\rho_{N+1}}\left(V_{\zeta_i,\zeta_{i+1}}V_{\zeta_j,\zeta_{j+1}}\right)\,,
$$
i.e. the correlations for the measure $\mu_N$ coincides with the correlations of the local function $V$ for the Markov measure $\rho_{N+1}$. Other features of the measure $\mu_N$ can be deduced by those of $\rho_{N+1}$ as for example large deviations \cite{GI} and FKG type inequalities \cite{CFGGT}.
\end{remark}

\begin{remark}\label{brno}
Item $(2)$ in the proof of Theorem \ref{ilteo} is a statement similar to the result in \cite{FNS}. In an infinite stationary framework, the authors of \cite{FNS} introduce the class of algebraic measures and the class of manifestly positive algebraic measures, for which probabilities of cylinder sets can be computed, like in \eqref{defbase}. They prove an equivalence between the family of manifestly positive algebraic measures and the family of functions of Markov processes. A process $(Y_n)_{n\in \mathbb N}$ is called a function of a Markov process if there exists a Markov process  $(X_n)_{n\in \mathbb N}$ and a function $\Phi$ such that $Y_n=\Phi(X_n)$. Also the authors of \cite{FNS} for simplicity consider finite dimensional cases and obtain bounds on the size of the configuration space of the Markov process. In the case that the manifestly positive algebraic measure has state space $A$ and the probabilities are computed using matrices of size $|B|$ then the process can be obtained as a function of a Markov process whose configuration space is bounded by $\left(\max\{|A|,|B|\}\right)^4$.

By item $(2)$ of Theorem \ref{ilteo} we can deduce that the process $(\eta_i)_{i=1}^N$ is a function of the inhomogeneous Markov process $\big(\eta_i,\zeta_i\big)_{i=1}^N$ via the function $\Phi(a,b)=a$. In this case the size of the configuration space of the Markov process is given by $|A||B|$; however by using the remaining statements of Theorem \ref{ilteo} we can instead obtain the mixture representation where several features of the process $(\eta_i)_{i=1}^N$ can be deduced from those of the Markov process with a state space whose cardinality is just $|B|$. This fact is simplifying very much the probabilistic description.
\end{remark}

\begin{remark}
	In the case of finite matrices, i.e. $|B|<+\infty$, the maximal positive eigenvalue coincides with the spectral radius of the matrix and there is only one positive eigenvector. This means that the transition matrix associated to a nonnegative matrix by \eqref{Doob} is unique. In the case of infinite matrices the spectral radius does not coincide with the maximal eigenvalue but has a different definition, see for example \cite{K}. In this case there may be more than one positive eigenvalue, all of them are associated to real positive eigenvalues bigger than the spectral radius; you may also have no eigenvectors. There is unicity in the case that the stochastic matrix constructed by \eqref{Doob} by using as eigenvalue the spectral radius is recurrent. We do not enter into details, see for example \cite{K} for some statements. In the infinite case we may therefore have more than one stochastic matrix associated by \eqref{Doob} to $M$ and this can be related to Doob $h$ transformations and also to existence of positive non constant harmonic functions. Suppose for example we constructed $P$ from $M$ by \eqref{Doob} and that $h$ is a non constant harmonic function for $P$ (i.e. satisfies $Ph=h$); then we have that $P^h$  defined by $P^h_{b,b'}=h_b^{-1}P_{b,b'}h_{b'}$, is related to $M$ by a generalized Doob transform like \eqref{Doob}. By the same telescopic arguments used in the proof of Theorem \ref{ilteo} we have that the families of Markov bridges constructed by using $P$ and $P^h$ coincide.
\end{remark}

\begin{remark}
	Finally we point out that Markov bridges are a remarkable class of processes for which a large literature is available, for example concerning also simulations.
\end{remark}

\section{examples}
In this section we apply the above constructions to specific examples considering also countable alphabets. The constructions are the same of the finite case, we need only to verify summability conditions that indeed are satisfied in the specific examples considered. Other examples will be discussed in \cite{I}.

\subsection{Boundary driven TASEP}

We consider the boundary driven totally asymmetric exclusion process (TASEP) on a one dimensional chain $\{1,2\dots ,N\}$ with N sites. This is a continuous time Markov process on $\{0,1\}^N$. An element of the state space $\eta=(\eta_1, \dots ,\eta_N)\in \{0,1\}^N$ represents a configuration of particles and is so that $\eta_i=1$ if there is a particle at site $i$ and $\eta_i=0$ if the site $i$ is empty. The generator is given by
\begin{align}\label{gen-ex}
\mathcal L_N f(\eta)&=\sum_{i=1}^{N-1}\eta_i(1-\eta_{i+1})\left[f(\eta^{i,i+1})-f(\eta)\right]\\
&+\alpha\left[f(\eta^{1,+})-f(\eta)\right]+\beta\left[f(\eta^{N,-})-f(\eta)\right]\,,
\end{align}
where  $0\leq\alpha,\beta\leq 1$ are parameters,
\begin{equation}\label{exchange}
\eta^{i,j}_k=\left\{
\begin{array}{ll}
\eta_k & \textrm{if}\ k\neq i,j \\
\eta_i & \textrm{if}\ k=j \\
\eta_j & \textrm{if}\ k=i
\end{array}
\right.
\end{equation}
and
\begin{equation}
\eta^{i,+}_k=\left\{
\begin{array}{ll}
\eta_k & \textrm{if}\ k\neq i\\
1 & \textrm{if}\ k=i\,
\end{array}
\right.
\qquad
\eta^{i,-}_k=\left\{
\begin{array}{ll}
\eta_k & \textrm{if}\ k\neq i\\
0 & \textrm{if}\ k=i\,.
\end{array}
\right.
\end{equation}

Informally we have a chain of $N$ sites where on each site at most one particle can be present. Particles jump with rate one to the nearest neighbor site to the right (i.e. from site $i$ to site $i+1$) and the jump is suppressed if this site is already occupied. Finally particles are injected with rate $\alpha$ on the leftmost site $1$, when empty, and particles in the rightmost site $N$ are destroyed with rate $\beta$.

This is an important non equilibrium model whose invariant measure can be represented in MPA form \eqref{defbase} with $A=\{0,1\}$ and $B=\mathbb N\cup\{0\}:=\mathbb N_0$, see \cite{Der3, Schutz}. Formula \eqref{defbase} is the invariant measure of this model, provided the following algebraic relations are satisfied by the two matrices $M^0, M^1$ for some vectors $x,y$
\begin{equation}\label{comm}
\left\{
\begin{array}{l}
M^1M^0=M^1+M^0 \,,\\
M^1 x=\frac 1\beta x\,, \\
y^T M^0=\frac 1\alpha y^T\,.
\end{array}
\right.
\end{equation}
Equations \eqref{comm} do not identify uniquely the matrices $(M^a)_{a=0,1}$ and there are several different solutions of \eqref{comm} for different values of the parameters.

We now apply Theorem \ref{ilteo} to some representations of the solutions of \eqref{comm}.

\subsubsection{The case $\alpha+\beta>1$}\label{secTASEP1}
In this case relations \eqref{comm} can be satisfied considering $B=\mathbb N_0$ and by the following matrices
\begin{equation}\label{matrici}
M^0=\begin{bmatrix}
1 & 0 & 0 & \cdots & 0 & \cdots \\ 1 & 1 & 0 & \cdots & 0 & \cdots \\ 0&1&1&0&\cdots &\cdots \\ 0& 0 &1&1&0 &\cdots \\ \vdots&\vdots&0&\ddots&\ddots&\ddots \\ \vdots&\vdots&\vdots&\ddots&\ddots&\ddots
\end{bmatrix}\text{,} \ \	M^1=\begin{bmatrix}
1 & 1 & 0 & \cdots & 0 & \cdots \\ 0 & 1 & 1 & 0& \cdots & \cdots \\ 0&0&1&1& 0& \cdots & \\ 0& 0 &0&1&1 &\ddots \\ \vdots&\vdots&\vdots&0&\ddots&\ddots \\ \vdots&\vdots&\vdots&\ddots&\ddots&\ddots
\end{bmatrix}
\end{equation}
and the vectors $y=\Big(\left(\frac{1-\alpha}{\alpha}\right)^b\Big)_{b\in \mathbb N_0}$, $x=\Big(\left(\frac{1-\beta}{\beta}\right)^b\Big)_{b\in \mathbb N_0}$, see\cite{Der3}.
The matrix $M=M^0+M^1$ is
\begin{equation}\label{M}
M_{b,b'}=\left\{
\begin{array}{ll}
2\delta_{b,b'}+\delta_{b,b'+1}+\delta_{b,b'-1}\,, & b\geq 1\,,\\
2\delta_{b,b'}+\delta_{b,b'+1}\,, & b=0\,,
\end{array}
\right.
\end{equation}
and the normalization factor \eqref{norm} can be easily computed and is finite so that the measure \eqref{defbase} is well defined .

\begin{figure}
\begin{center}
	\includegraphics[scale=1.8]{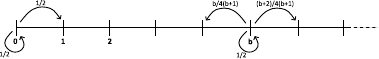}
\end{center}
\caption{The transition graph associated to the matrix $P$ in \eqref{bT1}.}\label{figura1}
\end{figure}

\begin{figure}
	\begin{center}
		\includegraphics{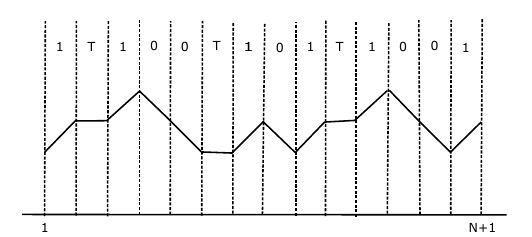}
	\end{center}
	\caption{A trajectory $\zeta$ of a Markov bridge for a random walk on $\mathbb N_0$ with jumps $0,\pm 1$. Given the trajectory, the configuration $\eta$ is completely determined on the intervals where the slope is $\pm 1$ while we need to toss independent fair coins (T=Toss) in correspondence of horizontal steps. The trajectory of $\zeta$ may also assume the value zero.}\label{figT}
\end{figure}

The eigenvalue problem for an eigenvalue $\lambda$ with corresponding eigenvector $(e_b)_{b\in \mathbb N_0}$ is given by
\begin{equation}\label{Tl}
\left\{
\begin{array}{ll}
2e_0+e_1=\lambda e_0\,, & \\
e_{i-1}+2e_i+e_{i+1}=\lambda e_i & i\geq 1\,,
\end{array}
\right.
\end{equation}
that is satisfied by $\lambda=4$ with the corresponding positive eigenvector $e=(b+1)_{b\in \mathbb N_0}$. Indeed it is possible to show that $\lambda=4$ is the spectral radius of the matrix $M$ (see \cite{K} for definitions and more), but in order to implement the transformation \eqref{Doob} we need just the algebraic relation \eqref{Tl}. The stochastic matrix \eqref{Doob} can therefore be explicitly computed and it is given by
\begin{equation}\label{bT1}
P_{b,b'}=\left\{
\begin{array}{ll}
\frac 12 & b'=b\\
\frac{(b+2)}{4(b+1)} & b'=b+1\\
\frac{b}{4(b+1)} & b'=b-1\,,
\end{array}
\right.
\end{equation}
and its transition graph is represented in Figure \ref{figura1}
We have also
\begin{equation}\label{fgTAS}
\left\{
\begin{array}{ll}
f_b=\frac{x_b}{b+1}=\frac{(1-\beta)^b}{\beta^b(b+1)} & b\in \mathbb N_0\,,\\
g_b=y_b(b+1)=\frac{(1-\alpha)^b(b+1)}{\alpha^b} & b\in \mathbb N_0\,.
\end{array}
\right.
\end{equation}
Finally, since $A=\{0,1\}$ and recalling that $\mathcal B_q\in \mathcal M^1(\{0,1\})$ is the Bernoulli measure of parameter $q$, using \eqref{defp} we have
\begin{equation}\label{lepT}
p_{b,b'}=\left\{
\begin{array}{ll}
\mathcal B_1 & \textrm{if}\ b'=b+1\,,\\
\mathcal B_0 & \textrm{if}\ b'=b-1\,,\\
\mathcal B_{\frac 12} & \textrm{if} \ b'=b\,.
\end{array}
\right.
\end{equation}
This completes the probabilistic characterization of the invariant measure of boundary driven TASEP for this values of the parameters $\alpha, \beta$. The Markov bridge $\rho_{N+1}$ for the variables $\zeta$ is associated to the transition matrix \eqref{bT1} and the vectors \eqref{fgTAS}, while the marginal laws of the mixtures are given by \eqref{lepT}. The construction of a sample of $\eta$ distributed according to the invariant measure starting from a sample of $\zeta$ is illustrated in Figure \ref{figT}.

\smallskip

\paragraph{\bf Comparison with the two lines ensemble}

Let us now show how the mixture representation that we obtained is related and equivalent to the  representation recently derived in \cite{Bryc} based on the approach of \cite{BC}. This is obtained using a pair of right-up independent random walks $s^1=(s^1_i)_{i=1}^{N+1}\,,\, s^2=(s^2_i)_{i=1}^{N+1}$ and is therefore called a two lines ensemble. Both representations describe the same invariant measure and must be of course equivalent in this sense. The two are however based on an enlargement of the state space with the introduction of a coupling with some hidden variables; what we show here is that the coupled measures, that could in principle be different, are instead the same modulo a simple mapping. There are several other combinatorial representation of the invariant measure of TASEP, see for example \cite{V, CW1, CW2, DS} and references therein.

\smallskip

The two lines ensemble is as follows.
We have that $s^1_1=s^2_1=0$ and for $2\leq i\leq N+1$ we have $s^1_i=\sum_{j=1}^{i-1}\eta_j$, $s^2_i=\sum_{j=1}^{i-1}\gamma_j$, where $\eta,\epsilon \in \{0,1\}^N$. We called $\eta$ the increments of the walk $s^1$ since they will be distributed finally according to the invariant measure of TASEP. The two lines ensembles representation in \cite{Bryc,BC} is written in terms of the two function $s^1,s^2$ defining the coupling measure
\begin{equation}\label{Brycf}
C_N^{\textrm{2-lines}}\left(s^1,s^2\right):=\frac{\left(\frac{1-\beta}{\beta}\right)^{\left\{s^1_{N+1}-s^2_{N+1}\right\}}}{\mathcal Z_N\left(\frac{(1-\beta)(1-\alpha)}{\beta\alpha}\right)^{\min_{1\leq i\leq N+1}\left\{s^1_{i}-s^2_{i}\right\}}}\left(\frac{1}{4}\right)^N\,,
\end{equation}
where $\mathcal Z_N$ is a normalization factor and $\left(1/4\right)^N=\mathbb P_{\textrm{Unif.}}(s^1,s^2)$ is
the uniform measure on the pair of paths that we do not include in the normalization factor just for convenience.
In \cite{Bryc} it is proved that $\mu_N^{\textrm{2-lines}}(s^1)=\sum_{s^2}C_N^{\textrm{2-lines}}\left(s^1,s^2\right)$, coincides with $\mu_N(\eta)$, by using the bijection between $s^1$ and $\eta$.

\smallskip

Let us introduce the $\mathbb N_0\times\mathbb N_0$ stochastic matrix $\hat P$ defined by
\begin{equation}\label{bT1hat}
\hat P_{b,b'}=\left\{
\begin{array}{ll}
\frac 12 & b'=b\\
\frac{1}{4} & b'=b+1\\
\frac{1}{4} & b'=b-1\,, b>0\,,\\
\frac{1}{2} & b=0, b'=1\,.
\end{array}
\right.
\end{equation}
This stochastic matrix is related to \eqref{bT1} by the relations
$$
\hat P_{b,b'}=\left\{
\begin{array}{ll}
\frac{P_{b,b'}e_{b}}{e_{b'}}\,, & b>0\,,\\
P_{b,b'}\,, & b=0\,.
\end{array}
\right.
$$
This means that the stochastic matrix $\hat P$ is almost related by a Doob transformation to the matrix $P$; the only transition on which the relation fails is the transition $0\to 1$. Consider $\zeta\in \mathbb N_0^{N+1}$ and call $\mathcal N_{0,1}(\zeta):=\sum_{i=1}^{N}\id \big(\zeta_i^{i+1}=(0,1)\big)$. By the above relations we have
\begin{equation}
\mathbb P_P(\zeta)=\mathbb P_{\hat P}(\zeta)\frac{\zeta_{N+1}+1}{\zeta_{1}+1}\left(\frac 12\right)^{\mathcal N_{0,1}(\zeta)}\,, \qquad \zeta\in \mathbb N_0^{N+1}\,.
\end{equation}
This allows to write the Markov bridge for the stochastic matrix $P$ in terms of the stochastic matrix $\hat P$ as
\begin{equation}\label{nuovaforma}
\rho_{N+1}(\zeta)=\frac{\left(\frac{1-\alpha}{\alpha}\right)^{\zeta_1}\mathbb P_{\hat P}(\zeta)
\left(\frac{1-\beta}{\beta}\right)^{\zeta_{N+1}}\left(\frac 12\right)^{\mathcal N_{0,1}(\zeta)}}{g^TP^Nf}\,,
\end{equation}
where the denominator is simply the normalization factor.

Consider a random walk $\xi_i\in \mathbb Z$ that jumps to $\xi_i\pm 1$ with probability $1/4$ and stays at $\xi_i$ with probability $1/2$. Its $\mathbb Z\times \mathbb Z$ transition probability $\tilde P$ is given by
\begin{equation}\label{bT1tilde}
\tilde P_{b,b'}=\left\{
\begin{array}{ll}
\frac 12 & b,b'\in \mathbb Z\,,\ b'=b\\
\frac{1}{4} & b,b'\in \mathbb Z\,, \ b'=b\pm 1\,.
\end{array}
\right.
\end{equation}
The stochastic matrix $\hat P$ is the transition probability of the Markov process $\zeta_i=|\xi_i|$. We have
in addition that for any $\zeta\in \mathbb N_0^{N+1}$
\begin{equation}\label{invtr}
\mathbb P_{\hat P}(\zeta)\left(\frac 12\right)^{\mathcal N_{0,1}(\zeta)}=
\mathbb P_{\tilde P}(b+\zeta) \qquad \forall b\in \mathbb Z\,,
\end{equation}
where by $b+\zeta$ we denote the shifted path $\big(b+\zeta_i\big)_{i=1}^{N+1}$ (this means that $\mathbb P_{\tilde P}(b+\zeta)=\left[\prod_{i=1}^N \tilde P_{b+\zeta_i,b+\zeta_{i+1}}\right]$).
The arbitrary shift factor follows by the fact that $\tilde P$ is associated to a space homogeneous random walk on $\mathbb Z$.

Recall the paths $(s^1_i)_{i=1}^{N+1}$ and $(s^2_i)_{i=1}^{N+1}$ introduced for defining the two lines ensembles and call $\xi_i=s^1_i-s^2_i$. We can describe the independent uniform measure on the two paths both in terms of the pairs $(s^1,s^2)$ or of the pairs $(\eta,\xi)$, where we recall the $\eta$ are the increments of $s^1$ and $\xi$ has just been defined. With a small abuse of notation we call $\mathbb P_{\textrm{Unif.}}(\eta,\xi)$ the measure $\mathbb P_{\textrm{Unif.}}(s^1,s^2)$ when written in terms of the variables $(\eta,\xi)$. By a direct computation we have
\begin{equation}\label{uqua}
\mathbb P_{\textrm{Unif.}}(\eta,\xi)=\mathbb P_{\tilde P}(\xi)\mathbb P_{\textrm{Unif.}}(\eta|\xi)=\mathbb P_{\tilde P}(\xi)\prod_{i=1}^Np_{\xi_i,\xi_{i+1}}(\eta_i)\,,
\end{equation}
where the probability measures $p_{b,b'}(\cdot)$ are those defined in \eqref{lepT}. An important features of such measures is that they are invariant by a joint shift of the indices, this means that $p_{b,b'}(\cdot)=p_{b+c,b'+c}(\cdot)$ for any $c\in \mathbb Z$. This means that the Markov bridge \eqref{nuovaforma} or any translation of such a measure are good parameter laws of the mixture in order to obtain the invariant measure of boundary driven TASEP.

Fix $\xi=(\xi_i)_{i=1}^{N+1}\in \mathbb Z^{N+1}$ such that $\xi_1=0$ and $|\xi_i-\xi_{i+1}|\leq 1$ and compute the weight given to all the translated of $\xi$ by the probability measure \eqref{nuovaforma}. We have therefore to sum the probability given by \eqref{nuovaforma} to all the paths of the form $\zeta=b+\xi$ where $b+\xi$ as before is the path shifted by $b\in \mathbb Z$. Since $\zeta\in \mathbb N_0^{N+1}$ we have the constraint that $b\geq -\min_{1\leq i\leq N+1}\xi_i:=-m_{N+1}(\xi)$. We obtain the probability measure $\tilde \rho_{N+1}$ on paths on $\mathbb Z$ starting from the origin defined by
\begin{eqnarray}
&\tilde\rho_{N+1}(\xi):=\sum_{\left\{b\geq -m_{N+1}(\xi)\right\} }\rho_{N+1}(b+\xi)\\
&=\frac{\mathbb P_{\tilde P}(\xi)\sum_{\left\{b\geq -m_{N+1}(\xi)\right\} }\left(\frac{1-\alpha}{\alpha}\right)^{b}
	\left(\frac{1-\beta}{\beta}\right)^{b+\xi_{N+1}}}{g^TP^Nf}\\
&=\frac{\mathbb P_{\tilde P}(\xi)
	\left(\frac{1-\beta}{\beta}\right)^{\xi_{N+1}}}{\left(\frac{(1-\alpha)(1-\beta)}{\alpha\beta}\right)^{m_{N+1}(\xi)}(g^TP^Nf)\left(1-\frac{(1-\alpha)(1-\beta)}{\alpha\beta}\right)}\,.
\end{eqnarray}
where we used \eqref{invtr} in the second equality while the last equality follows by a direct computation of the geometric sum.
By the translational invariance of \eqref{lepT} we have that the coupled measure
\begin{equation}\label{newC}
C_N(\eta,\xi)=\tilde\rho_{N+1}(\xi)\prod_{i=1}^Np_{\xi_i,\xi_{i+1}}(\eta_i)
\end{equation}
is so that $\mu_N(\eta)=\sum_\xi C_N(\eta,\xi)$ is the invariant measure of boundary driven TASEP. Using \eqref{uqua} and writing the measure in terms of the paths $(s^1,s^2)$, we obtain the two lines ensemble \eqref{Brycf} with $\mathcal Z_N=\left(g^TP^Nf\right)\left(\frac{\alpha+\beta-1}{\alpha\beta}\right)$.

\subsubsection{A general case}
\begin{figure}
	\begin{center}
		\includegraphics[scale=1.8]{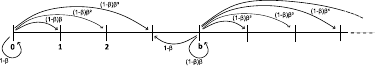}
	\end{center}
	\caption{The transition graph on $\mathbb N_0$ for the transition matrix \eqref{Pconb}; the Markov dynamics may have arbitrary long jumps.}\label{fig2}
\end{figure}

A different solution of the algebraic relations \eqref{comm} that holds for any value of the parameters, is obtained again with $B=\mathbb N_0$ and by the matrices
\begin{equation}
M^0=\begin{bmatrix}
0 & 0 & 0 & \cdots & 0 & \cdots \\ 1 & 0 & 0 & \cdots & 0 & \cdots \\ 0&1&0&0&\cdots &\cdots \\ 0& 0 &1&0&0 &\cdots \\ \vdots&\vdots&0&\ddots&\ddots&\ddots \\ \vdots&\vdots&\vdots&\ddots&\ddots&\ddots
\end{bmatrix}\text{,} \ \	M^1=\begin{bmatrix}
\frac{1}{\beta} & \frac{1}{\beta} & \frac{1}{\beta} & \cdots & \frac{1}{\beta} & \cdots \\ 0 & 1 & 1 & 1& \cdots & \cdots \\ 0&0&1&1& 1& \cdots & \\ 0& 0 &0&1&1 &\ddots \\ \vdots&\vdots&\vdots&0&\ddots&\ddots \\ \vdots&\vdots&\vdots&\ddots&\ddots&\ddots
\end{bmatrix}
\end{equation}
with vectors $y=\Big(\frac{1}{\alpha^b},\Big)_{b\in \mathbb N_0}$, and $x=\Big(\delta_{b,0}\Big)_{b\in \mathbb N_0}$, see again \cite{Der3}. The matrix $M=M^0+M^1$ is:
\begin{equation}
M_{b,b'}=\left\{
\begin{array}{ll}
\frac{1}{\beta} &b=0\,, \\
1\, & \forall \ b' \geq b-1\,.
\end{array}
\right.
\end{equation}
The eigenvalue problem for an eigenvalue $\lambda$ and a corresponding eigenvector $e=(e_b)_{b\in \mathbb N_0}$ is obtained by solving the equations
\begin{equation}
\left\{
\begin{array}{ll}
\lambda e_0=\frac 1\beta \sum_{j=0}^{+\infty}e_j\,, & \\
\lambda e_i=\sum_{j=i-1}^{+\infty} e_j\,, & i\geq 1\,.
\end{array}
\right.
\end{equation}
As it is possible to check a positive eigenvalue is given by $\lambda=\frac{1}{\beta(1-\beta)}$ with corresponding positive eigenvector given by $e=(\beta^b)_{b \in \mathbb N_0}$. We can then apply the transformation \eqref{Doob} obtaining explicitly the corresponding stochastic matrix that is given by:
\begin{equation}\label{Pconb}
P_{b,b+k}=\left\{
\begin{array}{ll}
(1-\beta)\beta^{k} &b=0, \ k \geq 0\,, \\
(1-\beta)\beta^{k+1} & \forall \ b \geq 1, \ k \geq -1\,,
\end{array}
\right.
\end{equation}
and all the remaining entrances are zero. The transition
graph of the stochastic matrix $P$ on $\mathbb N_0$ is drawn in Figure \ref{fig2}. The vectors defining the Markov bridge are $g=\Big(\left(\frac{\beta}{\alpha}\right)^b\Big)_{b\in \mathbb N_0}$ and $f=\big(\delta_{0,b}\big)_{b\in \mathbb N_0}$.
The marginal distributions \eqref{defp} can also be easily computed and are given by:
\begin{equation}\label{Plong}
p_{b,b'}=\left\{
\begin{array}{ll}
\mathcal B_0 &b'=b-1\,, \\
\mathcal B_1 & \forall \ b' \geq b\,.
\end{array}
\right.
\end{equation}
Recall that $\mathcal B_q$ is a Bernoulli measure of parameter $q$, and since the parameters are only $1$ and $0$, this means that the particle configuration $\eta$ is completely determined by the Markov bridge $\zeta$.

\smallskip

\paragraph{\bf An exploration process}

We discuss now a natural and simple interpretation of the Markov process obtained above.
There is a well known exploration algorithm of branching processes (see for example section 3.3 in \cite{vdh}) that connect branching processes to random walks on $\mathbb N_0$ that may have arbitrary large positive jumps but negative jumps only equal to $-1$. The Markov bridge obtained above is constructed starting from Markov chains of this type. We consider a small variation of the classic construction; in particular we have an exploration of infinite independent branching processes each of them attached to its own root vertex. Attach to any root an independent copy of a branching process with branching law $q\in \mathcal M^1(\mathbb N_0)$. Consider now the random walk $(\zeta_n)_{n\in \mathbb N}$ associated to the following exploration procedure. Start with $\zeta_0=0$ and select one root. The value $\zeta_1$ is the number of nodes attached to the root and set all of them in the {\it activated} status. Select one node among the activated ones, according for example to the depth-first search and set it in the {\it explored} status. This node does not belong any more to the set of activated nodes (it has been explored) and add to the set of activated nodes all the vertices in the progeny of the node just explored. The variable $\zeta_2$ is the number of activated nodes after this construction. We can now iterate the procedure and call $\zeta_n$ the number of activated nodes after $n$ iterations. If after n steps all nodes of the tree associated to the branching process have been explored then $\zeta_n=0$ since there are no more activated nodes. You proceed then in the next step to the exploration of the branching tree attached to another root and the variable $\zeta_{n+1}$ is the number of nodes attached to the new root, and so on. The random walk that we described has the transition probability equal to $P$ in \eqref{Pconb} when the branching law $q$ is a geometric of parameter $(1-\beta)$ that means $q_k=(1-\beta)\beta^k$, $k=0,1,2..$

\begin{figure}
	\begin{center}
		\includegraphics[scale=1.1]{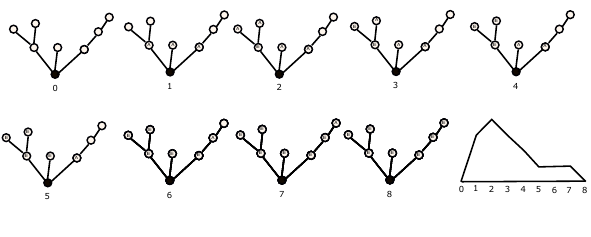}
	\end{center}
	\caption{The exploration sequence of a rooted tree, the root is the black node. The symbol $A$ means Activated while the symbol $E$ means Explored; at the end it is drawn the path of the corresponding walk $\zeta$ associated to the exploration procedure.}\label{figex}
\end{figure}

\subsection{TASEP with second class particles on a ring}

We consider the TASEP on $\mathbb Z_N$ the  ring with $N$ sites where particles can be of two different type, first and second class particles. We have that $A=\{0,1,2\}$ and the value of $\eta_i$ represents the class of the particle  at $i$; this means that if $\eta_i=0$ then the site $i\in \mathbb Z_N$ is empty, if $\eta_i=1$ there is a first class particle at $i$, while when $\eta_i=2$  there is a second class particle at $i$. Since we are on a ring we use an equivalence modulo $N$ for the indices on the lattice.
The stochastic dynamics is defined by the generator
\begin{equation}
\mathcal L_Nf(\eta):=\sum_{i=1}^N \id(r(\eta_i)>r(\eta_{i+1}))\left[f(\eta^{i,i+1})-f(\eta)\right]\,,
\end{equation}
where the symbol $\eta^{i,i+1}$ is defined in \eqref{exchange} and $r:\{0,1,2\}\to \{0,1,2\}$ is the priority function defined by: $r(0)=0\,, r(1)=2\,, r(2)=1$; particles of higher class have lower priority. The informal description of the dynamics is as follows: at rate one the occupation variables at the extreme of a bond $(i,i+1)$ exchange their value when the priority of the variable at vertex $i$ is higher than that at vertex $i+1$.

The dynamics has two conserved quantities that are the number of first class particles (i.e. the number of particles with priority $2$) $\mathcal N_1(\eta):=\sum_{i=1}^N\id(\eta_i=1)$ and the number of second class particles (i.e. the number of particles with priority $1$) $\mathcal N_2(\eta):=\sum_{i=1}^N\id(\eta_i=2)$. This means that there is a two parameter family of invariant measures depending on the pair of integer numbers $(n_1,n_2)$ that represent respectively the number of first and second class particles and are such that $n_1+n_2\leq N$. We call $\mu_N^{n_1,n_2}$ the corresponding invariant measure when it is necessary to indicate the parameters.
The matrix product formulation on a periodic lattice is slightly different from \eqref{defbase} and we have that the stationary measure is given by \cite{Der2cl}:
\begin{equation}\label{mes2class}
\mu_N^{n_1,n_2}(\eta)=\frac{\text{Tr} \left[ \prod_{i=1}^{N} M^{\eta_i}\right]}{Z_N(n_1,n_2)}\id\big( \mathcal N_1(\eta)=n_1\,,\mathcal N_2(\eta)=n_2\big)\,,
\end{equation}
where the normalization factor $Z_N(n_1,n_2)$ depends on the numbers $n_1,n_2$ of the first and second class particles. The matrices $M^i$, $i=0,1,2,$ have to satisfy the following conditions:
\begin{equation*}
M^1M^0=M^1+M^0, \ \ \ \ \ M^1M^2=M^2, \ \ \ \ M^2M^0=M^2.
\end{equation*}
A solution to the above algebraic relations is obtained (see \cite{Der2cl}) with $B=\mathbb N_0$ and the matrices given by:
\begin{equation*}
M^0=\begin{bmatrix}
1 & 0 & 0 & \cdots & 0 & \cdots \\ 1 & 1 & 0 & \cdots & 0 & \cdots \\ 0&1&1&0&\cdots &\cdots \\ 0& 0 &1&1&0 &\cdots \\ \vdots&\vdots&0&\ddots&\ddots&\ddots \\ \vdots&\vdots&\vdots&\ddots&\ddots&\ddots
\end{bmatrix}\text{,} \ \	M^1=\begin{bmatrix}
1 & 1 & 0 & \cdots & 0 & \cdots \\ 0 & 1 & 1 & 0& \cdots & \cdots \\ 0&0&1&1& 0& \cdots & \\ 0& 0 &0&1&1 &\ddots \\ \vdots&\vdots&\vdots&0&\ddots&\ddots \\ \vdots&\vdots&\vdots&\ddots&\ddots&\ddots
\end{bmatrix}
\end{equation*}
\begin{equation*}
M^2=\begin{bmatrix}
1 & 0 & 0 & \cdots & \cdots \\ 0 & 0 & 0 & \cdots & \cdots \\ 0&0&0&0&\cdots \\ \vdots & \vdots &0&0&\cdots \\ \vdots&\vdots&\ddots&\ddots&\ddots \\
\end{bmatrix}
\end{equation*}
When $n_2=0$ and there are only first class particles, the measure \eqref{mes2class} is not defined, since the trace is infinite, but the invariant measure is easily computed directly as the uniform measure on all the compatible configurations. We assume therefore that $n_2>0$ that means that there is at least one second class particle and this is always understood when the indices $n_i$ appear. A probabilistic description of such measures has been given in \cite{FFK}.

\begin{figure}
	\begin{center}
		\includegraphics[scale=1.8]{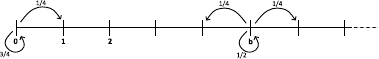}
	\end{center}
	\caption{The transition graph associated to the Markov matrix $P$ in \eqref{P2c}.}\label{fig3}
\end{figure}

We define the matrix $M=M^0+M^1+M^2$ that has the form
\begin{equation*}
M=\begin{bmatrix}
3 & 1 & 0 & \cdots & 0 & \cdots  & \cdots & \cdots\\ 1 & 2 & 1 & 0 & 0 & \cdots & \cdots & \cdots \\ 0&1&2&1&0 &\cdots & \cdots & \cdots \\ 0& 0 &1&2&1 &0 & \cdots & \cdots\\ \vdots&\vdots&0&1&2&1 & 0 & \cdots \\ \vdots&\vdots&\vdots&0 &\ddots&\ddots & \ddots & \ddots
\end{bmatrix}
\end{equation*}
The sum of each row of $M$ is constant and equal to $4$ so that $\lambda=4$ is a positive eigenvalue with corresponding constant eigenvector $e=(1)_{b\in \mathbb N_0}$. By the transformation \eqref{Doob} we obtain the stochastic matrix $P$ defined by
\begin{equation}\label{P2c}
P_{b,b'}=\left\{
\begin{array}{ll}
1/4 & b>0\,,\, b'=b\pm 1\,, \\
1/2 & b>0\,,\, b'=b\,, \\
1/4 & b=0\,,\, b'=1\,, \\
3/4 & b=0\,,\, b'=0\,,
\end{array}
\right.
\end{equation}
whose transition graph is represented in Figure \ref{fig3}.

The set of invariant probability measures is a convex set whose extremal elements are the measures in
\eqref{mes2class}. In \cite{Der2cl} the authors introduce a grand canonical invariant probability measure
\begin{equation}\label{mes2classN}
\mu_N(\eta)=
\frac{\text{Tr} \left[ \prod_{i=1}^{N} M^{\eta_i}\right]}{Z_N}\id\big(\mathcal N_2(\eta)>0\big)\,,
\end{equation}
that is a non trivial convex combination of the extremal invariant measures. We have indeed that $\mu_N=\sum_{n_1,n_2}\frac{Z_N(n_1,n_2)}{Z_N}\mu^{n_1,n_2}_N$, where the sum is over the pairs $(n_1,n_2)$ such that $n_1+n_2\leq N$ and $n_2>0$. The extremal measures can be obtained as the canonical measures of the grand canonical one
$$
\mu_N\left(\cdot |\mathcal N_1=n_1,\mathcal N_2=n_2\right)=\mu^{n_1,n_2}_N(\cdot)\,.
$$

\smallskip

To give a mixture representation of the invariant measures we need to introduce another probability measure, since the construction of Theorem \ref{ilteo} using the measure \eqref{mes2classN} does not work properly. This is due to the fact that if we start from the measure \eqref{mes2classN} and enlarge it to the coupling measure $C_N$, as we did before, our general strategy fails. In particular it is not possible to implement the marginal computation of the issue $3)$ of the proof of Theorem \eqref{ilteo}, since the coupled measure contains some constraints written in terms of the variables $\eta$.  This constraint is contained in the characteristic function and  as a result we have that the marginal $\zeta$ law of the coupling is not a simple Markov process. We construct therefore another grand canonical measure whose canonical measures are again the measures $\mu_N^{n_1,n_2}$.

Fix $j$ a node of the ring and let us define
\begin{equation}\label{atmuj}
\mu_{N,j}:=\frac{\text{Tr} \left[ \prod_{i=1}^{N} M^{\eta_i}\right]}{Z^*_N}\id\big(\eta_j=2\big)\,,
\end{equation}
where by symmetry the normalization factor $Z^*_N$ does not depend on $j$. We define also the translational invariant measure $\hat\mu_N:=\frac 1N \sum_{j=1}^N\mu_{N,j}$ for which we have
\begin{equation}\label{atmu}
\hat\mu_N(\eta)=\frac{\text{Tr} \left[ \prod_{i=1}^{N} M^{\eta_i}\right]\sum_{j=1}^N\id(\eta_j=2)}{NZ^*_N}=
\frac{Z_N\mu_N(\eta)\mathcal N_2(\eta)}{NZ^*_N}\,.
\end{equation}
By a direct computation it is easy to see that the measure \eqref{atmu} can be written as
$\hat \mu_N=\sum_{n_1,n_2}\frac{n_1Z_N(n_1,n_2)}{NZ^*_N}\mu^{n_1,n_2}_N$ (we recall once again that the sum is over $(n_1,n_2)$ such that $n_1+n_2\leq N$ and $n_2>0$) and in particular again \eqref{mes2class} are the conditional measures of $\hat\mu_N$.

The coupling construction of Theorem \ref{ilteo} is particularly simple for the measures $\mu_{N,j}$ and consequently for $\hat\mu_N$. The motivation is that the constraint on the $\eta$ variables can be transformed in this case into a constraint on the $\zeta$ variables. By symmetry the construction for the different $j$ is the same up to a rotation, for simplicity of notation we discuss the case $\mu_{N,N}$. For each $j$ we can construct the coupling $C_{N,j}\in \mathcal M^1(A^N\times B^N)$ that in the case $j=N$ is
\begin{align*}
C_{N,N}(\eta,\zeta)&:= \frac{\left[\prod_{i=1}^{N} M^{\eta_i}_{\zeta_i,\zeta_{i+1}}\right]\id\big(\eta_N=2\big)}{Z^*_N}\\
&=\frac{3\left[\prod_{i=1}^{N-1} M^{\eta_i}_{\zeta_i,\zeta_{i+1}}\right]\id\big(\eta_N=2,\zeta_1=0, \zeta_N=0\big)}{Z^*_N}\,.
\end{align*}
The $\zeta$ marginal of the coupled measure $C_{N,N}$ can be easily computed since the sum over the $\eta_i$ with $i\neq N$ has no constraints while there is no sum over $\eta_N$ that is fixed equal to $2$; we obtain
\begin{equation}\label{bridgeper}
\rho_{N,N}(\zeta)=\sum_\eta C_{N,N}(\eta,\zeta)=\frac{3\prod_{i=1}^{N-1} P_{\zeta_i,\zeta_{i+1}}}{4^{N-1}Z^*_N}\id\big(\zeta_1=0, \zeta_N=0\big)\,.
\end{equation}
In the case of periodic boundary conditions, since we are on a ring, the number of $\zeta$ variables is equal to the number of the $\eta$ variables and we have therefore a Markov bridge of length $N$ and not $N+1$ as in the open boundary case. According to the notation \eqref{defbr}, we have that \eqref{bridgeper} is the  Markov Bridge $\rho^{0,0}_N$ of length $N$ with transition matrix $P$ and pinned to start and finish at $\zeta_1=\zeta_N=0$. It is geometrically natural to associate the variables $\zeta$ to the dual lattice of that to which the variables $\eta$ are associated and in such a way that $\zeta_i$ and $\zeta_{i+1}$ are associated to the points of the dual lattice corresponding to the edges exiting from the vertex of the original lattice associated to the variable $\eta_i$. We will use a bit this shifted association in the discussion of the comparison with other representations. A similar construction can be done in the open boundary case. By rotational symmetry the measures $\rho_{N,j}$ have the same structure, after a rotation.

Also in this case the marginals of the product measures can be directly computed. For any $b,b'\in B$ we have $p_{b,b'}(\cdot)\in \mathcal M^1(\{0,1,2\})$ and we use the notation $p_{b,b'}=(p_{b,b'}(0),p_{b,b'}(1),p_{b,b'}(2))$ By computing \eqref{defp}  we obtain four different probability measures depending on the values of $b,b'$
\begin{equation}\label{p2cl}
\left\{
\begin{array}{ll}
p_{0,0}=\left(\frac 13, \frac 13, \frac 13\right)\,, & \\
p_{b,b+1}=(0,1,0)\,, & \\
p_{b,b-1}=(1,0,0)\,, & \\
p_{b,b}=\left(\frac 12,\frac 12,0\right)\,, & b\geq 1\,.
\end{array}
\right.
\end{equation}
To construct a configuration $\eta$ distributed according to the grand canonical invariant measure \eqref{atmu}, we put a second class particle at $N$, i.e. we fix $\eta_N=2$. Then we construct a sample path $\zeta$ according to the Markov Bridge $\rho_{N,N}=\rho_N^{0,0}$. Given $\zeta$ we generate the $\eta_i$, $i\neq N$, independently using \eqref{p2cl}. Once obtained the configuration $\eta$ in this way, we insert in the ring simply shifting by a uniform rotation.

\smallskip

\paragraph{\bf Comparison with the queue representation}

We start recalling the collapsing representation of the invariant measure \cite{A} and its queue interpretation and generalization in \cite{FM}. An interesting issue would be to consider the mixture representation in the multiclass case of \cite{FM} which representation in MPA form is obtained in \cite{EFM}.

\smallskip

We use the language of \cite{FM} with a left-right symmetry exchange of the argument due to the different orientation of the dynamics there. Consider two subsets $\mathcal A, \mathcal S \subseteq  \mathbb Z_N$ such that $|\mathcal A|< |\mathcal S|$; they may have elements in common. The set $\mathcal A$ is the set of arrivals of clients while the subset $\mathcal S$ is the set of services.  The time axis for the queue interpretation is going from right to left according to the same order of going from $i+1$ to $i$. We introduce the variables $(A_i)_{i\in \mathbb Z_N}$ and $(S_i)_{i\in \mathbb Z_N}$ defined by $A_i=1$ if $i\in \mathcal A$ and zero otherwise, while $S_i=1$ if $i\in \mathcal S$ and zero otherwise. The construction in \cite{A,FM} is as follows. We construct iteratively a subset $\mathcal S^+\subseteq \mathcal S$ that is the set of {\it used} service times. By construction we will have $|\mathcal S^+|=|\mathcal A|$ and we call $\mathcal S^-:=\mathcal S \setminus \mathcal S^+$ the set of {\it unused} services. Initially $\mathcal S^+$ is empty. Select an arbitrary client (an $i\in \mathcal A$) and serve it in the nearest unused service to its left (a $j\in \mathcal S\setminus \mathcal S^+$). The selected service time is added to the set $\mathcal S^+\subseteq \mathcal S$ of used services. You iterate the procedure up to assign a service time to any client. Given $\mathcal S^\pm$ we assign a configuration of particles putting a first class particle on each vertex in $\mathcal S^+$, a second class particle on each vertex in $\mathcal S^-$ and leaving empty the remaining vertices. When $|\mathcal S|=n_1+n_2$, $|\mathcal A|=n_1$ and they are uniformly distributed among all the sets with such constraints, then the configuration of particles obtained is distributed like $\mu_N^{n_1,n_2}$ in \eqref{mes2class}.

\smallskip

Since our mixture representation is associated to the measures \eqref{atmuj}, we need to slightly modify this construction in order to have a direct comparation, we consider the case $j=N$. We call $\mathbb Z^*_N$ the dual lattice. Consider $\mathcal S, \mathcal A$ two subsets of $\mathbb Z_N$ such that $|\mathcal A|<|\mathcal S|$.
We represent elements of $\mathbb Z_N$ periodically on $\mathbb Z$ and fix an arbitrary reference starting point $j\in \mathbb Z^*$ on the dual lattice. Define on $\mathbb Z^*$ the function $\tilde W_i$, $i\in \mathbb Z^*$, defined by
$$
\tilde W_i:=\left\{
\begin{array}{ll}
\sum_{\ell\in \mathbb Z\cap [i,j]}\left(A_\ell-S_\ell\right)\,, & i\leq j\,, \\
-\sum_{\ell\in \mathbb Z\cap [j,i]}\left(A_\ell-S_\ell\right)\,, & i> j\,.
\end{array}
\right.
$$
\begin{figure}
	\begin{center}
		\includegraphics[scale=0.8]{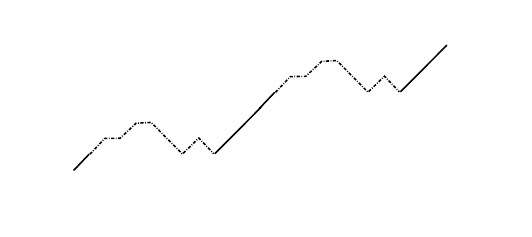}
	\end{center}
	\caption{An example of a function $\tilde W$ with the subdivision of the graph into excursions (dotted part) and records (continuous part).}\label{figrecord}
\end{figure}
By periodicity we have $\tilde W_{i+N}=\tilde W_i+|\mathcal S|-|\mathcal A|$ and then we have $\lim_{i\to \pm\infty}\tilde W_i=\pm\infty$, since $|\mathcal S|>|\mathcal A|$. We extend the function $\tilde W$ to a function defined on the whole real axis by linear interpolation. The graph of $\tilde W$ is divided into records and excursions; a piece of the graph between the coordinates $i\in \mathbb Z^*$ and $i+1\in \mathbb Z^*$ is a record if $\tilde W_i <\tilde W_k$ for any $k>i$ in $\mathbb Z^*$, otherwise belongs to an excursion, see Figure \ref{figrecord} for an illustrative example. By periodicity this classification is invariant by shift of $N$ and we obtain a division of the edges of $\mathbb Z^*_N$, that are in bijection with the nodes of $\mathbb Z_N$, into edges belonging to excursions and to records. The values of $(A_i,S_i)$ with $i\in \mathbb Z_N$ corresponding to an edge of the dual graph that is a record can be only $(0,1)$. Since $|\mathcal A|<|\mathcal S|$ and for each excursion we have the same number of values $(0,1)$ (that correspond to decreasing steps going from right to left in $\tilde W$) and $(1,0)$ (that correspond to increasing steps going from right to left in $\tilde W$), we obtain that there must exists at least one $j\in \mathbb Z_N$ such that $(A_j,S_j)=(0,1)$ and moreover $j$ is in correspondence to a record. We fix the coordinates on the torus in such a way that $j=N$. We consider the function $\tilde W$ constructed using as a starting reference point on the dual lattice exactly $j-\frac 12$ that changing the reference coordinates becomes $N-\frac 12$. We consider this function $\tilde W$ in the interval $\left[\frac 12,N-\frac 12\right]$ and define $W:=\left[\tilde W\right]_+$, where $[\cdot]_+$ denotes the positive part. By the special choice of the starting point
the function $W$ is simply obtained from $\tilde W$ by preserving the shape of excursions and transforming into constant pieces with value zero the records; in particular $W_{\frac 12}=W_{N-\frac 12}=0$. The function $W$ represent the length of the queue associated to the arrival and services determined respectively by $\mathcal A, \mathcal S$ and starting the counting from the initial reference point. Apart the $1/2$ shift due to a simpler and clearer description, the graph of $W$ has the same geometric features of the sample paths of $\zeta$ distributed according to $\rho_N^{0,0}$ and indeed we are going to show an identification.

\begin{figure}
	\begin{center}
		\includegraphics[scale=1.3]{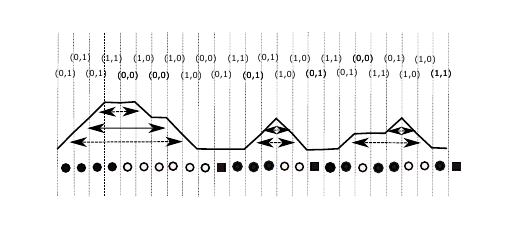}
	\end{center}
	\caption{An example of the identification of used and unused services that generates a configuration of particles of first and second class. In the top part of the picture there are the values of the variables $(A_i,S_i)$ and below them there is the graph of $W$ associated to the waiting queue for the service. The horizontal arrows show the matching among clients and services according to the last arrived first served procedure. Below we draw the configuration of particles obtained: round black dots are first class particles, black squares are second class particles and round white dots are empty sites. Note that there is an extra second class particle on the right of all, that corresponds to the coordinate $N$.}\label{figpairing}
\end{figure}

The identification of the sets $\mathcal S^\pm$ with the special reference frame that we selected is particularly simple. The length of the queue starting from the reference point is $W$. If at $i$ the queue is not zero, i.e. $W_i>0$, and there is a service time available, i.e. $S_i=1$, then we serve one client in the queue at this time; this means that $i$ is an used service time and $i\in \mathcal S^+$. Also if the queue at $i$ is zero but $(A_i,S_i)=(1,1)$ we have at $i$ an used service and $i\in \mathcal S^+$ (the client arrives and is immediately served). Note that any $i$ such that $(A_i,S_i)=(1,1)$ is in an excursion. All the remaining $i\in \mathcal S$ are instead unused services and belong to $\mathcal S^-$. This construction corresponds to identify the vertices in $\mathcal S^+$ with the elements of $\mathcal S$ that belong to excursions while the vertices in $\mathcal S^-$ with the elements of $\mathcal S$ that belongs to records. This is obtained just considering the unfair procedure of last arrived first served and we obtain the pairing inside each excursion illustrated in figure \ref{figpairing}.
\smallskip

We construct now starting from sets $\mathcal A$ and $\mathcal S$ the variables $\zeta$. We start constructing the unconditioned Markov chain with transition matrix \eqref{P2c} and then the variables $\eta$ that conditioned on the $\zeta$ have a product distribution with marginals \eqref{p2cl}. To compare directly with the construction in the general part we identify, like before, the element of the dual lattice between the vertices of the original lattice in $i$ and $i+1$ with $i$. Consider $A_i$ and $S_i$ that are independent i.i.d and uniformly distributed assuming values $0,1$ with the same probability.

We consider the walks $(s^1_i)_{i=1}^N$ defined by $s^1_1=0$ and $s^1_i=\sum_{\ell=1}^{i-1}S_i$, $i>1$; $(s^2_i)_{i=1}^N$ defined by $s^2_1=0$ and $s^2_i=\sum_{\ell=1}^{i-1}A_i$, $i>1$. We define then $\xi_i=s^1_i-s^2_i$ and finally $\zeta_i:=|\xi_i|$. We have that the variables $\zeta$ constructed in this way form a Markov chain with initial condition $\zeta_1=0$ and transition matrix \eqref{P2c}.
The variables $\eta$ can now be defined as follows in such a way that conditioned on the $\zeta$ they have law \eqref{p2cl}
\begin{equation}\label{basta}
\eta_i=S_{i+1}\id(\zeta_i>0)+\id(\zeta_i=0)\Big[2S_{i+1}-S_{i+1}A_{i+1}\Big]\,, \qquad i=1, \dots , N-1\,.
\end{equation}
Considering subsets $\mathcal S, \mathcal A$ for which $\zeta_N=0$ we see that formula \eqref{basta} and the queue construction done using the same subsets identify exactly the same splitting $\mathcal S=\mathcal S^+\cup \mathcal S^-$. This means that when $\mathcal A, \mathcal S$ have the correct conditional distribution the mixture construction and the queue one identify the same configuration of particles. We do not discuss more details, the equivalence is even more simple and direct in the infinite setting of $\mathbb Z$.

\subsection{Harmonic models}
These models have been introduced in \cite{FGK} and a combinatorial representation of the invariant measures has been obtained in \cite{FG}. We consider only the simplest case.
We have a one dimensional chain $\{1,\dots, N\}$ with two external sources at the left and right boundaries having two parameters $0<\alpha <\beta<1$, respectively (the case $\beta<\alpha$ can be considered just exchanging left with right, in the case $\alpha=\beta$ the model is reversible and the invariant measure is product). On each lattice site we have an arbitrarily large number of particles and we denote by $\eta_i\in \mathbb N_0$ the number (possibly zero) of particles at $i\in \{1,\dots, N\}$. We consider a continuous-time Markov chain $\{\eta(t), \ t\ge 0\}$ whose state space is the set $\mathbb N_0^N$ of configurations $\eta=(\eta_1,\dots ,\eta_N)$. The stochastic dynamics has a bulk and a boundary part which are described in terms of the generator $L_N$ defined below. For any $i\in \{1,\dots, N\}$, we denote by $\delta^i$ the configuration defined by $\delta^i_j=0$ when $j\neq i$ and $\delta^i_i=1$.
We have
\begin{equation}\label{generator}
L_N:=L_N^{\textrm{bulk}}+L_N^{\textrm{bound}}.
\end{equation}
The bulk generator applied to bounded functions reads:
\begin{equation}
\label{genbulk}
L_N^{\textrm{bulk}}f(\eta)= \sum_{\substack{i,j \\ |i-j|=1}} \sum_{k=1}^{\eta_i}\frac 1k \left[f(\eta-k\delta^i+k\delta^j)-f(\eta)\right]\,.
\end{equation}
The boundary part which encodes the interaction with the reservoirs
is given by:
\begin{align}
\label{genbound}
L_N^{\textrm{bound}}f(\eta) & = \sum_{k=1}^{\eta_1}\frac 1k \left[f(\eta-k\delta^1)-f(\eta)\right]+
\sum_{k=1}^{\infty}\frac{\alpha^k}{k}\left[f(\eta+k\delta^1)-f(\eta)\right]\nonumber \\
& + \sum_{k=1}^{\eta_N}\frac 1k \left[f(\eta-k\delta^{N})-f(\eta)\right]+
\sum_{k=1}^{\infty}\frac{\beta^k}{k}\left[f(\eta+k\delta^{N})-f(\eta)\right]\,.
\end{align}
Let $\mathcal G_m(k)=\frac{1}{1+m}\left(\frac{m}{1+m}\right)^k$, $k=0,1,\dots$, be a geometric distribution of mean $m$. Given $\underline m=(m_1,\dots ,m_{N})$ and $\underline k =(k_1, \dots ,k_{N})$ we denote by
$\mathcal G_{\underline m}(\underline k):=\prod_{i=1}^{N}\mathcal G_{m_i}(k_i)$.
Given $0<\alpha<\beta <1$ we call $m_L:=\frac{\alpha}{1-\alpha} < \frac{\beta}{1-\beta}:=m_R$ and introduce
$O^{m_L,m_R}_N\subseteq [m_L,m_R]^{N}$ as the set defined by
$$
O^{m_L,m_R}_N:=\left\{\underline m\,:\,m_L\leq m_1\leq \dots \leq m_{N}\leq m_R\right\}\,,
$$
whose Lebesgue volume is given by $|O^{m_L,m_R}_N|=\frac{(m_L-m_R)^{N}}{N!}$.
Inspired by \cite{BDGJL,CGT}, in \cite{CFGGT} it has been proved that the invariant measure of this model can be represented as
	\begin{equation}\label{formulasuper}
	\mu_{N}^{m_L,m_R}(\eta)=\frac{1}{|O^{m_L,m_R}_N|}\int_{O^{m_L,m_R}_N} d \underline m\  \mathcal G_{\underline m}(\eta)\,.
	\end{equation}
For convenience we make explicit the dependence of the invariant measure on the parameters $m_L, m_R, N$.

\smallskip

Formula \eqref{formulasuper} represents the invariant measures as a mixture of product of inhomogeneous geometric distributions. We show now that formula \eqref{formulasuper} can be interpreted as a probability measure of the MPA type \eqref{defbase} but with operators $(M^k)_{k\in \mathbb N_0}$ that are determined by continuous kernels $\left(M^k_{m,m'}\right)_{m,m'\in \mathbb R^+}^{k\in \mathbb N_0}$. Using our equivalence we then show that \eqref{formulasuper} is a particular case of \eqref{mixform} for a special Markov bridge and special marginal distributions \eqref{spform}. Recently a representation of \eqref{formulasuper} as a measure of MPA type with matrices has been obtained in \cite{Fr}; it would be interesting to explore the corresponding mixture representation that should be different from \eqref{formulasuper} (recall Remark \ref{Remo1} on non uniqueness).

\smallskip

Let us define the kernels
\begin{equation}\label{operator}
M^k_{m,m'}:=\mathcal G_m(k)\id (m\leq m')\,,\qquad  k\in \mathbb N_0\,,\  m,m'\in \mathbb R_+\,,
\end{equation}
that are kernels of operators in a functional space that for simplicity we do not describe formally. We use the bra-ket formalism of quantum mechanics to denote elements of the space on which the operators act. For a state $\ket{f}$ associated to a function $f:\mathbb R_+\to \mathbb R_+$ we have
\begin{equation}
\left\{
\begin{array}{l}
M^k\ket f=\int_0^{+\infty}dm'\, M^k_{m,m'}f(m')=\mathcal G_m(k)\int_m^{+\infty}f(m')d\,m'\,,\\
\bra f M^k=\int_0^{+\infty}dm\, M^k_{m,m'}f(m)=\int_0^{m'}\mathcal G_m(k)f(m) d\,m\,.
\end{array}
\right.
\end{equation}
We have also delta states $\ket z$ associated to elements $z\in \mathbb R_+$ and the action of the operators is defined in this case as
\begin{equation}
\left\{
\begin{array}{l}
M^k\ket z=\int_0^{+\infty}dm'\, M^k_{m,m'}\delta(m'-z)=\mathcal G_m(k)\id (m\leq z)\,,\\
\bra z M^k=\int_0^{+\infty}dm\, M^k_{m,m'}\delta(m-z)=\mathcal G_z(k)\id (z\leq m')\,.
\end{array}
\right.
\end{equation}
Using these operators we can introduce the family of probability measures
\begin{equation}\label{muh}
\hat\mu_{N+1}^{m_L,m_R}(\eta):=\frac{\bra{m_L} \prod_{i=0}^N M^{\eta_i}\ket{m_R}}{Z_N}\,,
\end{equation}
where the states $\bra m_L$ and $\ket m_R$ are delta states associated to the values $m_L, m_R \in \mathbb R_+$. The index $N+1$ is due to the fact that in \eqref{muh} we have $\eta=(\eta_0,\eta_1,\dots ,\eta_N)$ containing the extra variable $\eta_0$ that is not belonging to our particle system;
we have however the relation
$$\mu^{m_L,m_R}_N(\eta_1,\dots ,\eta_N)=\sum_{\eta_0}\hat\mu_{N+1}^{m_L,m_R}(\eta_0,\dots ,\eta_N)\,.$$
The operator $M=\sum_{k\in \mathbb N_0} M^k$ has a continuous family of positive eigenvectors; consider the family of functions $\left(e_\lambda(m)=\lambda e^{-\lambda m}\,; m\in \mathbb R_+\right)_{\lambda >0}$, then we have
$$
Me_\lambda=\int_m^{+\infty}\lambda e^{-\lambda m'} \, dm'=e^{-\lambda m}= \frac{e_\lambda}{\lambda}\,.
$$
For any $\lambda>0$ we can therefore apply the transformation \eqref{Doob} obtaining the transition kernel
\begin{equation}\label{ph}
P^\lambda_{m,m'}=\lambda e^{-\lambda (m'-m)}\id (m'-m\geq 0)\,d\,m'\,.
\end{equation}

Since the states in \eqref{muh} are delta states, the Markov bridge $(\zeta_0,\dots ,\zeta_{N+1})$ is a Markov bridge with transition probabilities given in \eqref{ph} and pinned at initial and final time at $\zeta_0=m_L$ and $\zeta_{N+1}=m_R$. The Markov bridge has length $N+2$ since in \eqref{muh} there is an extra variable $\eta_0$. We have therefore that the random variables $(\zeta_0,\dots, \zeta_{N+1})$ are distributed as the random walk $\zeta_n=\zeta_0+\sum_{i=1}^n \gamma_i$, where the $(\gamma_i)_{i\in \mathbb N}$ are i.i.d. exponential random variables of parameter $\lambda>0$. For different values of $\lambda$ we have different Markov processes but having the same conditioned Markov bridge, that indeed coincides with the order statistics law of uniform on the interval $[m_L,m_R]$. A special feature of these Markov processes is that they are increasing, so that the Markov property is  equivalent to the spatial Markov property of the point process $\{\zeta_1,\dots , \zeta_N\}\subseteq [m_L,m_R]$ discussed in \cite{CFF,RV}; this is not true in the general case.

By a simple direct computation we obtain that the marginals \eqref{defp} of the mixture are given by
$p_{m,m'}(k)=\mathcal G_m(k)$; also in this case this model has again a special feature that simplify the computation and the representation, that is the fact that the marginal of the mixture depends just on $m$ and not on $m'$.

A similar generalized construction as above can be done for the whole class of harmonic models whose representation of the invariant measures as a mixture has been obtained in \cite{CFF}.

\subsection{Miscellany}

In the finite dimensional case, under irreducibility conditions, there is a unique stochastic matrix $P$ associated to $M$ by \eqref{Doob}. This means that, not considering for the moment the vectors $x,y$, the whole class of probability measures of MPA type \eqref{defbase}, can be parametrized by a stochastic matrix $P$ and $|B|^2$ elements of $\mathcal M^1(A)$ that determine the marginals $p_{b,b'}(\cdot)$, $b,b'\in B$ of the mixture. This is a smaller family with respect to the one obtained by the family of $|A|$ non negative  $B\times B$ matrices $M^a$. The mixture representation \eqref{mixform} can be therefore interpreted as a canonical form of the measures of MPA type; given $P$ and $\left(p_{b,b'}\right)_{b,b'\in B}$ there is a whole family of matrices $(M^a)_{a\in A}$ that are associated to them and are all those of the form
\begin{equation}\label{famiglia}
M^a_{b,b'}=p_{b.b'}\lambda \,e_b P_{b,b'}e_{b'}^{-1}\,, \lambda\in \mathbb R_+\,, e\in \mathbb R_+^B\,,\qquad a\in A\,,\, b,b'\in B\,,
\end{equation}
where $\lambda\in \mathbb R_+$ and $e\in \mathbb R_+^B$ (that can be fixed in such a way that $\sum_b e_b=1$) are arbitrary parameters.

Relation \eqref{famiglia} identifies equivalence classes of family of matrices $(M^a)_{a\in A}$ all of them corresponding to the same mixture. Note that this is not equivalent to the problem of the different representations of the matrices satisfying for example the algebraic relations \eqref{comm}. In that case we have infinite dimensional representations and moreover the matrices solving \eqref{comm} are not invariant by an arbitrary scaling factor, as instead is the case of \eqref{famiglia}.

The Markov bridge measure is determined also by the measure $m$ in \eqref{emme} that depends on the vectors $f,g$ that are related to the vectors $x,y$ by the matrices $\mathcal E\,,\, \mathcal E^{-1}$. It can be shown that, for a given $P$, not all the probability measures can be written as \eqref{emme} but any pair of one dimensional marginals can indeed be obtained, for large enough $N$ and $|B|<+\infty$. The vectors $f, g$ can be fixed imposing the normalization $\sum_bf_b=\sum_b g_b=1$. We do not discuss further this topic.

\medskip

We say that $\mu_N$ of MPA type is stationary when $\mu_N(\eta_i)=\mu_N(\eta_j)$ for any $i,j$. In this case we can extend the measure to a measure on infinite sequences. Using the mixture representation this happens when the Markov bridge is stationary and this happens when $g_b=\pi_b$ and $f_b=1$, $b\in B$, where $\pi$ is the invariant measure of $P$. In the case that for example $|A|=|B|=2$ the $\zeta$ process is a stationary $\{0,1\}$ valued Markov chain and the $\eta$ process is constructed, using 4 numbers $0<p_{b,b'}<1$, starting from a sample path of $\zeta$ and tossing independent coins of parameter $p_{b,b'}$ in correspondence of occurrences of $bb'$ in $\zeta$. Among these processes there are the special cases of example 5 of section 3.1 of \cite{FNS}. Differently from \cite{FNS} where you need to suitably enlarge the alphabets to describe the processes as functions of Markov processes, the description in terms of mixture keep the original sizes of $A,B$.

\medskip

Always in the stationary case, we consider the {\it one dependent} measures of example 4 in section 3.1 of \cite{FNS}. This corresponds to measures for which $\eta_i^{j-1}$ and $\eta_{j+1}^k$ for any $i<j<k$ are independent, that corresponds to have $\sum_{\eta_j}\mu_N(\eta)=\mu_N(\eta_1^{j-1})\mu_N(\eta_{j+1}^N)$.

As explained in \cite{FNS} a family of matrices $(M^a)_{a\in A}$ and vectors $x,y$, gives a stationary one dependent measure when $M=cxy^T$ is a rank one matrix for an arbitrary constant $c>0$. In this case the Perron eigenvector is $x$ with corresponding eigenvalue $cy^Tx$ and the stochastic matrix \eqref{Doob} is $P_{b,b'}=\frac{y_{b'}x_{b'}}{y^Tx}$ that corresponds to an i.i.d. process.
We obtain therefore a simple characterization of one dependent measures that are also manifestly positive algebraic measures, in the terminology of \cite{FNS}. These are measures obtained by a double independent construction parametrized by a $q\in \mathcal M^1(B)$ and a family $p_{b,b'}\in \mathcal M^1(A)$, $b,b'\in B$. Before you construct a sample path of $\zeta$ using $i.i.d.$ variables distributed as $q$; conditioned on this sample path you generate independent variables $\eta$ such that $\eta_i$ is distributed as $p_{\zeta_i,\zeta_{i+1}}$. The result is a one dependent stationary measure.

\subsection*{Acknowledgements}
We thank P.A. Ferrari, M. Goldwurm, F. Mignosi, D. Schmid and T. Prosen for useful discussions. We thank T.Prosen for pointing out reference \cite{FNS}.

D.G. acknowledges the financial support from the Italian Research Funding Agency (MIUR) through
PRIN project ``Emergence of condensation-like phenomena in interacting particle systems: kinetic and lattice models'', grant n. 202277WX43.

F.I. acknowledges the financial support of the INdAM - GNFM Project with code CUP $\sharp E5324001950001\sharp $.

\end{document}